\documentclass[11pt]{article}

\usepackage{amssymb}
\usepackage{amsmath}
\usepackage{theorem}
\usepackage{epsfig}
\usepackage{verbatim}
\usepackage{graphicx}
\usepackage{indentfirst}
\newtheorem{thm}{Theorem}[section]

\newtheorem{prop}[thm]{Proposition}

\newtheorem{rem}[thm]{Remark}

\newcommand{\R}{\Bbb{R}}

\newcommand{\T}{\mathbb{T}}
\newcommand{\D}{\displaystyle}

\newcommand{\curl}{{\rm curl}\thinspace}

\newcommand{\grad}{\nabla}

\newcommand{\dxu}{\partial_{x_1}}
\newcommand{\dxd}{\partial_{x_2}}

\newcommand{\dpt}{\partial_t}
\newcommand{\dist}{{\rm dist}\thinspace}

\newcommand{\la}{\Lambda}

\newcommand{\al}{\alpha}
\newcommand{\ep}{\varepsilon}

\newcommand{\dpa}{\partial^{\bot}_{\alpha}}
\newcommand{\da}{\partial_{\alpha}}
\newcommand{\g}{\gamma}

\newcommand{\ff}{\mathbb}
\newcommand{\pa}{\partial}
\newcommand{\va}{\varpi}
\newcommand{\si}{\sigma}

\def\eop{{\ \vrule height 3pt width 3pt depth 0pt}}

\numberwithin{equation}{section}

\begin{document}

\author{Angel Castro, Diego C\'ordoba and Francisco Gancedo}

\title{A naive parametrization for the vortex-sheet problem.}
\date{}

\maketitle

\begin{abstract}
We consider the dynamics of a vortex sheet that evolves by the
Birkhoff-Rott equations. The fluid evolution is understood as a weak
solution of the incompressible Euler equations where the vorticity
is given by a delta function on a curve multiplied by an
amplitude.  The solutions we study are with finite energy, which implies
zero mean amplitude. In this context we choose a parametrization for the motion of the vortex sheet for which
the equation is well-posed for analytic initial data. For the equation of the
amplitude we show ill-posedness for non-analytic initial data.
\end{abstract}

\maketitle

\section{Introduction}

We consider a velocity field $v = (v_1,v_2)$ satisfying the incompressible 2D- Euler equations
\begin{equation}\label{euler}
v_t+(v\cdot\grad) v=-\nabla p,
\end{equation}
\begin{equation}
\label{div0}\nabla\cdot v = 0.
\end{equation}
We study weak solutions of the system whose vorticity $\omega = \nabla\times v$ is a delta function supported on the curve $z(\al,t)$:
\begin{equation}\label{fwintro}
\omega(x,t)=\varpi(\al,t)\delta(x-z(\al,t)),
\end{equation}
i.e. $\omega$ is a measure defined by
\begin{equation*}
<\omega,\eta>=\int \varpi(\al,t)\eta(z(\al,t)) d\al,
\end{equation*}
with $\eta(x)$ a test function.

Here we shall assume for the curve the following scenarios:
\begin{itemize}
\item Periodicity in the horizontal space variable: $z(\al+
2k\pi,t) = z(\al,t)+(2k\pi,0)$.
\item A closed contour: $z(\al+
2k\pi,t) = z(\al,t)$.
\item An open contour vanishing at infinity: $\D\lim_{\al\rightarrow\infty}(z(\al,t)-(\al,0))=0$.
\end{itemize}

The vortex sheet $z(\alpha,t)$ evolves satisfying the equation
\begin{align}
\begin{split}\label{fullpm}
z_t(\al,t)&=BR(z,\varpi)(\al,t)+c(\al,t)\da z(\al,t),
\end{split}
\end{align}
where the  Birkhoff-Rott integral on the curve, which  comes from Biot-Savart law, is given by
\begin{align}
\begin{split}\label{fibr}
BR(z,\varpi)(\al,t)=\frac{1}{2\pi}
PV\int \frac{(z(\al,t)-z(\beta,t))^{\bot}}{|z(\al,t)-
z(\beta,t)|^2}\varpi(\beta,t)d\beta,
\end{split}
\end{align}
and $c(\alpha,t)$ represents the re-parametrization freedom.
Then we can close the system using Bernoulli's law with the equation:
\begin{align}
\begin{split}\label{fw}
\varpi_t&=\da (c\, \varpi).
\end{split}
\end{align}

%

We study initial data for the amplitude of the vorticity with mean zero which is preserved by equation \eqref{fw}. From Biot-Savart law, at first expansion,  the expression at infinity is of the order of $\frac{1}{|x|}\int\va$ for a closed curve o near planar at infinity. To obtain a velocity field in $L^2$ it is necessary to have $\int\va=0$ (for more details see \cite{bertozzi-Majda}). In the periodic case, $z(\al+2\pi k,t)=z(\al,t)+(2\pi k,0)$,
the following classical identity for complex numbers$$
\frac{1}{\pi}\big(\frac1z+\sum_{k\geq 1}\frac{2z}{z^2-(2\pi
k)^2}\big)=\frac{1}{2\pi\tan(z/2)},
$$yields (ignoring the variable $t$)
\begin{align*}
v(x)\!=\!\frac{-1}{4\pi}
\int_{-\pi}^\pi\!\!\!\!\varpi(\beta)\big(\frac{\tanh(\frac{x_2-z_2(\beta)}{2})
(1\!+\!\tan^2(\frac{x_1-z_1(\beta)}{2}))}{\tan^2(\frac{x_1-
z_1(\beta)}{2})\!+\!\tanh^2(\frac{x_2-z_2(\beta)}{2})}     ,\!\frac{\tan(\frac{x_1-z_1(\beta)}{2})
(\tanh^2(\frac{x_2-z_2(\beta)}{2})\!-\!1)}{\tan^2(\frac{x_1-
z_1(\beta)}{2})\!+\!\tanh^2(\frac{x_2-z_2(\beta)}{2})}\big)d\beta,
\end{align*} for $x\neq z(\al,t)$. Then $$\lim_{x_2\rightarrow\pm\infty}v(x,t)=\mp\frac1{4\pi}\int_{-\pi}^\pi \va(\beta) d\beta (1,0),$$ and to have the same value at infinity it is necessary again mean zero.

The problem of existence of weak solutions of the Euler equations for general initial velocity in $L^2$ is not well understood \cite{bertozzi-Majda}. There is solution for this problem but the velocity field becomes a Laplace-Young measure (see \cite{DM2}). Constantin, E and Titi \cite{CET} prove a condition of regularity in 3D within the chain of Besov spaces, $v\in L^3([0,T];B^{\al,\infty}_3)\cap C([0,T];L^2)$ with $\al>1/3$, for  weak solutions conserving energy
(Onsager's conjecture). Nevertheless there are results of non-uniqueness for weaker solutions with zero initial data that becomes nontrivial (see \cite{Scheffer} and \cite{Shnirelman}) even for velocity fields in $L^2$, i.e.  $v(x,t)\in L^{\infty}_c([0,T]; L^2)$ (see \cite{Camilo}).  There  is also a result of uniqueness for a vorticity in $L^1\cap L^\infty$ due to Yudovich \cite{Y}.

For the particular case of a vortex sheet there are many papers which consider the case of $\va$ with a distinguished sign. We can point out the work of Delort \cite{Delort} where he proves global existence of weak solutions for initial velocity in $L^2_{loc}$ and vorticity a positive Radon measure. A simpler proof can be found in \cite{Majda} due to Majda. Existence for a particular case of a Radon measure with non distinguished sign is shown in \cite{lopes1} . Also, in the case of analytic initial data, a local existence result for the vortex sheet is given by Sulem, Sulem, Bardos and Frisch in \cite{SSBF} in the case where the curve is represented by a graph. The first result of ill-posedness in Sobolev spaces for amplitude with a distinguished sign is due to Ebin \cite{Ebin2} in a bounded domain. In the same year Duchon and Robert \cite{DR} proved global-existence for peculiar initial data. They consider a particular $c(\al,t)$ which gives $z_{1t}(\al,t)=0$ and therefore if one parametrizes initially $z_0(\al)=(\al,y_0(\al))$ the free boundary is given in terms of a function and the equations \eqref{fullpm} become
\begin{align*}
\begin{split}
y_t(\al,t)&=\frac{1}{2\pi}PV\int\frac{(\al-\beta)+(y(\al,t)-y(\beta,t))\da y(\al,t)}{(\al-\beta)^2+(y(\al,t)-y(\beta,t))^2}\va(\beta,t) d\beta\\
\end{split}
\end{align*}
and $c(\alpha,t)$ in equation (\ref{fw}) is given by
$$
c(\al,t)=\frac{1}{2\pi}PV\int\frac{(y(\al,t)-y(\beta,t))}{(\al-\beta)^2+(y(\al,t)-y(\beta,t))^2}\va(\beta,t) d\beta.
$$
A similar approach is done by Caflish and Orellana \cite{CO} to show also global-existence for particular initial data and moreover they give an argument to prove ill-posedness in
$H^s$ for $s>3/2$. They choose $c(\al,t)=0$ which implies $\va(\al,t)=\va_0(\al)$. If $\va_0(\al)$ has a distinguish sign, the following change of variable is legitimate
$$
d\sigma=\va_0(\beta)d\beta
$$
and equations \eqref{fullpm} can be written as
\begin{align}
\begin{split}\label{B-R}
z_t(\al,t)&=\frac{1}{2\pi}PV\int\frac{(z(\al,t)-z(\beta,t))^{\bot}}{|z(\al,t)-z(\beta,t)|^2} d\beta,
\end{split}
\end{align}
which is the Birkhoff-Rott equation. By taking $z(\al,t)=(\al+\ep_1(\al,t),\ep_2(\al,t))$ (or $\va(\al,t)=1+\ep_1(\al,t)$ and $y(\al,t)=\ep_2(\al,t)$ in the parametrization of Duchon and Robert) and linearizing in \eqref{B-R} one can obtain
$$\partial_t\ep_1=-\frac12\la(\ep_2),\qquad\partial_t\ep_2=-\frac12\la(\ep_1),$$
where $\Lambda$ is the operator $\Lambda=(-\Delta)^\frac{1}{2}$.
Therefore
$$\widehat{\ep}_1(\xi,t)=\frac{\widehat{\ep}_1(\xi,0)+\widehat{\ep}_2(\xi,0)}{2}e^{-\pi|\xi|t}+\frac{\widehat{\ep}_1(\xi,0)-\widehat{\ep}_2(\xi,0)}{2}e^{\pi|\xi|t},$$
$$\widehat{\ep}_2(\xi,t)=\frac{\widehat{\ep}_1(\xi,0)+\widehat{\ep}_2(\xi,0)}{2}e^{-\pi|\xi|t}-\frac{\widehat{\ep}_1(\xi,0)-\widehat{\ep}_2(\xi,0)}{2}e^{\pi|\xi|t}.$$
Since the initial data $\ep_1(\xi,0)=\ep_2(\xi,0)$ only oscillate the dissipative waves, it follows global-existence even for non-regular initial data. Applying Fourier techniques to the nonlinear case
$$\partial_t\ep_1=-\frac12\la(\ep_2)+T(\ep_1,\ep_2),\qquad\partial_t\ep_2=-\frac12\la(\ep_1)+S(\ep_1,\ep_2),$$
 yields that these particular initial data, small enough, activate only the dissipative waves and control the nonlinear operators $T$ and $S$ obtaining global in time solutions.

The main idea to show ill-posedness of Caflisch and Orellana is to consider the following function
$$s_0(\gamma,t)=\ep(1-i)[(1-e^{-t/2-i\gamma})^{1+\nu}-(1-e^{-t/2+i\gamma})^{1+\nu}]$$
which is a solution of the linearization of equation (\ref{B-R}). For $0<\nu<1$, $s_0$ has a infinite curvature at $\gamma=t=0$. Then they prove that a function $r(\gamma,t)$ exists  such that $z(\gamma,t)=\gamma+s_0+r$ is an analytic solution of equation (\ref{B-R}) with infinite curvature at $\gamma=t=0$. Then they obtain ill-posedness in Sobolev spaces in the Hadamard sense using the following symmetry properties:

If $z(\gamma,t)$ is a solution of (\ref{B-R}) then so are $z_b(\gamma,t)=\overline{z}(\gamma,-t)$, $z_s(\gamma,t)=z(\gamma,t-t_0)$ and $z_n(\gamma,t)=n^{-1}z(n\gamma,nt).$

A study of the existence of solutions of equation (\ref{B-R}) in less regular spaces than $H^s$ can be found in \cite{Wu}. We also quote that the first evidence of singularities with analytic initial data was given by Moore in \cite{Moore}.
$$
$$

In this paper our first step will be to deduce the equation of motion of the vortex sheet from the weak formulation of the Euler equations. In order to do that, in proposition \ref{presiones}, we show the equality of pressure when we take the limits at each side of the curve for weak solutions satisfying \eqref{fwintro} (see also \cite{Roman}). In sections 3 and 4 we shall study the case in which the term in the tangential direction is given by $c(\al,t)=\frac12H(\va)(\al,t)$, where $H\va$ is the Hilbert transform of the function $\va$ (see \cite{St1}) given by
$$H\va(\alpha)=\frac{1}{\pi}PV\int_{-\infty}^{\infty} \frac{\va(\beta)}{\alpha-\beta}d\beta,$$
and in the periodic domains also by
$$H\va(\alpha)=\frac{1}{\pi}PV\int_{-\pi}^\pi \frac{\va(\beta)}{2\tan((\al-\beta)/2)}d\beta.$$
This term is just of the same order that the Birkhoff-Rott integral over $\va$ for a regular one-to-one curve \cite{ADY}. In fact, for $z(\al,t)=(\al,0)$, we have exactly
$$BR(z,\va)(\al,t)=\frac12H(\va)(\al,t)(0,1).$$

In  section 3 we show that the chosen parametrization provides solutions of the vortex sheet problem. The requirements are the usual for this system: the initial data have to be analytic. With our analysis we do not need to parameterize the interface in terms of a function as in \cite{SSBF}, the initial curve has to be one-to-one and with nonzero tangent vector. In the argument we modify the proofs used in the Cauchy-Kowalewski
theorems given in \cite{Nirenberg} and \cite{Nishida} in order to deal with the arc-chord condition.

Finally, in section 4,  we show ill-posedness for the equation of the amplitude (\ref{fw}) for non-analytic initial data with mean zero.

\section{The evolution equation}

In this section we shall obtain the classic Birkhoff-Rott
equations from the weak formulation of the Euler equations for the
velocity.  For this purpose we use the continuity of the pressure
over the vortex sheet (see the appendix below).  One alternative
equation deduction, which does not need to prove this property of
the pressure, can be found in \cite{lopes2} where the authors use
the weak formulation of the Euler equations for the vorticity in
2D. We consider weak solutions of the system
(\ref{euler}--\ref{div0}): for any smooth functions $\eta$ and
$\zeta$ compactly supported on $[0,T)\times\R^2$, i.e. in the
space $C^\infty_c([0,T)\times\R^2)$, we have
\begin{equation}\label{weakeuler}
\int_0^T\int_{\R^2}\big(v\cdot (\eta_t + v\cdot\nabla\eta)
+p\nabla\cdot\eta\big)dxdt+\int_{\R^2}v_0(x)\cdot\eta(x,0)dx = 0
\end{equation}
and
\begin{equation}\label{weakdiv0}
\int_0^T\int_{\R^2} v\cdot \nabla \zeta dxdt= 0,
\end{equation}
where $v_0(x)= v(x,0)$ is the initial data. 
This formulation is easily seen to be equivalent to the
more usual one with div-free test fields, which
does not explicitly contain the pressure.


Let us assume that the vorticity is given by a delta function on the curve $z(\alpha,t)$ multiplied by an amplitude, i.e.
\begin{equation}\label{omegacomoes}
\omega(x,t)=\varpi(\al,t)\delta(x-z(\al,t)),
\end{equation}
where $z(\alpha,t)\in C^{1,\delta}$ splits the plane in two domains $\Omega^j(t)$ ($j=1,2$) and $\varpi(\al,t)\in C^{1,\delta}$ with $0<\delta<1$.

 Then by Biot-Savart law we get
\begin{equation}\label{vfc}
v(x,t)=\frac{1}{2\pi}PV\int\frac{(x-z(\beta,t))^\bot}{|x-z(\beta,t)|^2}
\varpi(\beta,t)d\beta
\end{equation}
for $x\neq z(\al,t)$. We have
\begin{align}\label{limitevelocidad}
\begin{split}
v^2(z(\al,t),t)&=BR(z,\varpi)(\al,t)+\frac12\frac{\varpi(\al,t)}{|\da
z(\al,t)|^2}\da z(\al,t),\\
v^1(z(\al,t),t)&=BR(z,\varpi)(\al,t)-\frac12\frac{\varpi(\al,t)}{|\da
z(\al,t)|^2}\da z(\al,t),
\end{split}
\end{align}
where $v^j(z(\alpha,t),t)$ denotes the limit velocity field obtained
approaching the boundary in the normal direction inside $\Omega^j$
and $BR(z,\varpi)(\al,t)$ is given by \eqref{fibr}. It is easy to
check that the velocity $v$ \eqref{vfc} satisfies \eqref{weakdiv0} .

 Next we shall obtain the equation for the curve $z(\alpha,t)$. We start from equation (\ref{weakeuler}) with $\eta(x,0)=0$, which is
\begin{equation}
\int_0^T\int_{\R^2}[v\cdot(\eta_t+v\cdot\nabla\eta)+p\nabla\cdot\eta]\,\, dx dt=0. \label{weakeuler2}
\end{equation}

Again we can split the equation (\ref{weakeuler2}) in the following way
$$\lim_{\ep \rightarrow 0}\Big( \int_0^T\int_{\Omega^1_{\ep}(t)}[v\cdot(\eta_t+v\cdot\nabla\eta)+p\nabla\cdot\eta]\,\, dx dt+$$
$$\int_0^T\int_{\Omega^2_{\ep}(t)}[v\cdot(\eta_t+v\cdot\nabla\eta)+p\nabla\cdot\eta]\,\, dx dt\Big)=0,$$
where
$$\Omega_{\ep}^1(t)=\{x\in\Omega^1(t):\dist(x,\partial\Omega^1(t))\geq\ep\}$$
$$\Omega_{\ep}^2(t)=\{x\in\Omega^2(t):\dist(x,\partial\Omega^2(t))\geq\ep\}.$$
We will study the first terms in detail.
Integrating by parts we obtain
$$\lim_{\ep\rightarrow 0}\int_0^T \int_{\Omega^1_{\ep}(t)}v\cdot\eta_t\,dxdt=\int_0^T\int (v^1\cdot\eta) (z_t\cdot \dpa z) d\alpha dt-\int_0^T\int_{\Omega^1(t)}v_t\cdot\eta\,dxdt.$$
Similarly for  the other terms we have
$$\lim_{\ep\rightarrow 0}\int_0^T\!\int_{\Omega^1_{\ep}(t)}v\cdot(v\cdot\nabla)\eta\,dxdt=-
\int_0^T\!\int (v^1\cdot\eta)(v^1\cdot \dpa z)\,d\alpha dt-\int_0^T\!\int_{\Omega^1(t)}\eta\cdot(v\cdot\nabla)v\,dxdt.$$
Operating in a similar way with the integral over $\Omega^2_{\ep}(t)$ yields the following equations
\begin{eqnarray}
\va\cdot(\dpt z-BR(z,\va))\cdot \dpa z=0,\\
v_t+(v\cdot\nabla)\cdot v=-\nabla p\quad \text{over $\Omega^1$ and $\Omega^2$},
\end{eqnarray}
where the derivatives of $v$ on $\partial\Omega^1$ and $\partial\Omega^2$ have to be understood like the  limits in the normal direction to the curve $z(\alpha,t)$ and we have used the continuity of the pressure (see Appendix).

Next we close the system giving the evolution equation for the amplitude of the
vorticity $\varpi(\al,t)$ by means of Bernoulli's law. Using
\eqref{vfc} for $x\neq z(\al,t)$ we get $v(x,t)=\grad \phi(x,t)$
where
$$
\phi(x,t)=\frac{1}{2\pi}PV\int
\arctan\Big(\frac{x_2-z_2(\beta,t)}{x_1-z_1(\beta,t)}\Big)\varpi(\beta,t)d\beta.
$$
We define
$$
\Pi(\al,t)=\phi^2(z(\al,t),t)-\phi^1(z(\al,t),t),
$$
where again $\phi^j(z(\alpha,t),t)$ denotes the limit obtained
approaching the boundary in the normal direction inside $\Omega^j$.
It is clear
\begin{align*}
\da \Pi(\al,t)&=(\grad\phi^2(z(\al,t),t)-\grad\phi^1(z(\al,t),t))\cdot \da z(\al,t)\\
&=(v^2(z(\al,t),t)-v^1(z(\al,t),t))\cdot \da z(\al,t)\\
&=\varpi(\al,t),
\end{align*}
therefore
$$
\int \varpi(\al,t) d\al=0.
$$
Now we can check that
\begin{align}\label{limitepotencial}
\begin{split}
\phi^2(z(\alpha,t),t)&=IT(z,\varpi)(\al,t)
+\frac{1}{2}\Pi(\al,t)\\
\phi^1(z(\alpha,t),t)&=IT(z,\varpi)(\al,t)-\frac{1}{2}\Pi(\al,t),
\end{split}
\end{align}
where
$$
IT(z,\varpi)(\al,t)=\frac{1}{2\pi}PV\int
\arctan\Big(\frac{z_2(\al,t)-z_2(\beta,t)}{z_1(\al,t)-z_1(\beta,t)}\Big)\varpi(\beta,t)d\beta.
$$

Using the Bernoulli's law in \eqref{euler}, inside each domain, we
have
$$
\phi_t(x,t)+\frac{1}{2}|v(x,t)|^2+p(x,t)=0.
$$
Taking the limit it follows
$$
\phi^j_t(z(\al,t),t)+\frac{1}{2}|v^j(z(\al,t),t)|^2+p^j(z(\al,t),t)=0,
$$
and since $ p^1(z(\al,t),t)=p^2(z(\al,t),t)$ we get
\begin{equation}\label{limiteB}
\phi^2_t(z(\al,t),t)-\phi^1_t(z(\al,t),t)+\frac{1}{2}|v^2(z(\al,t),t)|^2-\frac{1}{2}|v^1(z(\al,t),t)|^2=0.
\end{equation}
Then it is clear that
$\phi^j_t(z(\al,t),t)=\partial_t(\phi^j(z(\al,t),t))-z_t(\al,t)\cdot\grad\phi^j(z(\al,t),t)$
and using \eqref{limitevelocidad} together \eqref{limitepotencial} in
\eqref{limiteB} we obtain
\begin{align}
\begin{split}\label{eep}
\Pi_t(\al,t)&=\varpi(\al,t)(z_t(\al,t)-BR(z,\varpi)(\al,t))\cdot\frac{\da z(\al,t)}{|\da z(\al,t)|^2}.
\end{split}
\end{align}
Taking one derivative with respect to $\alpha$ on (\ref{eep}) yields equation (\ref{fw}).

Finally it is easy to show that the solutions of the system (\ref{fullpm}) and (\ref{fw}) provide weak solutions of the Euler's equation.

Given a curve $z(\alpha,t)\in C^{1,\,\delta}$ and a function $\va(\alpha,t)\in C^{1,\,\delta}$ such that the equations (\ref{fullpm}) and (\ref{fw}) are satisfied, we define the velocity $v(x,t)$ by the expression (\ref{vfc}) and the pressure by
$$p(x,t)=-\phi_t(x,t)-\frac{1}{2}|v(x,t)|^2\quad \text{over $\Omega^1$ and $\Omega^2$},$$
where the potential $\phi(x,t)$ is given by $v=\nabla \phi$.
From equation (\ref{fw}) we have that the pressure is continuous over the vortex sheet. In order to check that $v(x,t)$ and $p(x,t)$ are  weak solutions of Euler's equations we just have to introduce them in the first member of (\ref{weakeuler}) and (\ref{weakdiv0}) and integrate by parts.


\section{Local-existence for analytic initial data}


We have the evolution equation given by
\begin{align}
\begin{split}\label{ecca}
z_t&=BR(z,\va)+H\va\da z,\\
\va_t&=\da(\va H\va).
\end{split}
\end{align}
In this frame, we consider a scale of Banach spaces $\{X_r\}_{r\geq 0}$ given by periodic real functions that can be extended analytically on the complex strip
$B_{r}=\{\al+i\zeta: \al\in\T, |\zeta|<r\}$ with norm
$$
\|f\|_r=\max_{0\leq k\leq 2}\D\sup_{\al+i\zeta\in B_r}|\da^k f(\al+i\zeta)|_*+\D\sup_{\al+i\zeta\in B_r, \beta\in\T}\frac{|\da^2f(\al+i\zeta)-\da^2f(\al+i\zeta-\beta)|_*}{|\beta|^\delta},
$$
with $0<\delta<1$ and $|\cdot|_*$ the modulus of a complex number. We then obtain the following theorem.

\begin{thm}\label{analytic}
 Let $z^0(\al)$ be a curve satisfying the arc-chord condition
\begin{equation}\label{arc-chord}
\frac{|z^0(\al)-z^0(\al-\beta)|^2}{|\beta|^2}>\frac{1}{M^2},
\end{equation}
and $z^0(\al),\varpi^0(\al)\in X_{r_0}$ for some $r_0>0$. Then, there exist a time $T>0$ and $0<r<r_0$ so that there is a unique solution to \eqref{ecca} in $C([0,T];X_r)$.
\end{thm}

\begin{rem}
In the proof it is easy to check that the tangential term is not harmful to the evolution equation of the curve. In fact, it is the easier to deal with. Also, with solutions of the system \eqref{ecca}, by a reparametrization, one could recover solutions of the vortex sheets problem with the more usual choice of $c$ such as the one given by the lagrangian velocities or the one with $c=0$ (taking  $v=\frac{v^1+v^2}{2}$). A similar theorem follows for all these parametrizations.
\end{rem}

It is easy to check that $X_{r}\subset X_{r'}$ for $r'\leq r$ due to the fact that $\|f\|_{r'}\leq \|f\|_{r}$.
A simple aplication of the Cauchy formula gives
\begin{equation}\label{Cauchy}
\|\da f\|_{r'}\leq \frac{C}{r-r'}\|f\|_{r},
\end{equation}
for $r'<r$.

The equation \eqref{ecca} can be extended on $B_r$ as follows:
\begin{align}
\begin{split}
z_t(\al+i\zeta,t)&=F_1(z(\al+i\zeta,t),\varpi(\al+i\zeta,t)),\\
\varpi_t(\al+i\zeta,t)&=F_2(\varpi(\al+i\zeta,t)).
\end{split}
\end{align}
with
$$F_1(z,\varpi)=BR(z,\varpi)+H\varpi\da z,$$
and
$$F_2(\varpi)=\da (\varpi H\varpi).$$

\begin{prop}\label{Lip}
Consider $0\leq r'<r$ and the open set $O$ in $B_{r}$ given by
\begin{equation}\label{openO}
O=\{z,\varpi\in X_{r}: \|z\|_{r},\|\va\|_{r}< R, \D\inf_{\al+i\zeta\in B_{r}, \beta\in\T} G(z)(\al+i\zeta,\beta) >\frac{1}{R^2}\},
\end{equation}
with
\begin{equation}\label{arc-chord-A}
G(z)(\al+i\zeta,\beta)=\Big|\frac{(z_1(\al+i\zeta)-z_1(\al+i\zeta-\beta))^2+(z_2(\al+i\zeta)-z_2(\al+i\zeta-\beta))^2}{\beta^2}\Big|_*.
\end{equation}
Then the function $F=(F_1,F_2)$ for $F:O\rightarrow X_{r'}$ is a continuous mapping. In addition, there is a constant
$C_R$ (depending on $R$ only) such that
\begin{equation}\label{cota}
\|F(z,\varpi)\|_{r'}\leq \frac{C_R}{r-r'}\|(z,\varpi)\|_{r},
\end{equation}
\begin{equation}\label{casiL}
 \|F(z^2,\varpi^2)-F(z^1,\varpi^1)\|_{r'}\leq \frac{C_R}{r-r'}\|(z^2-z^1,\varpi^2-\varpi^1)\|_{r},
\end{equation}
and
\begin{equation}\label{paraarc-chord}
\sup_{\al+i\zeta\in B_r,\beta\in\T} |F_1(z,\varpi)(\al+i\zeta)-F_1(z,\varpi)(\al+i\zeta-\beta)|_*\leq C_R|\beta|,
\end{equation}
for $z,z^j,\varpi,\varpi^j\in O$.
\end{prop}

Using the above proposition we have the proof of theorem \ref{analytic}.\\

Proof of Theorem \ref{analytic}: The argument is analogous as in \cite{Nirenberg} and \cite{Nishida} (see also \cite{MP}). We have to deal with the arc-chord condition so we will point out the main differences. For initial data $z^0,\va^0\in X_{r_0}$ satisfying \eqref{arc-chord}, we can find a $0<r_0'<r_0$ and a constant $R_0$ such that
$\|z^0\|_{r_0'}< R_0$, $\|\va^0\|_{r_0'}< R_0$ and
\begin{equation}\label{arc-chord-}
\Big|\frac{(z^0_1(\al+i\zeta)-z^0_1(\al+i\zeta-\beta))^2+(z^0_2(\al+i\zeta)-z^0_2(\al+i\zeta-\beta))^2}{\beta^2}\Big|_*>\frac{1}{R_0^2},
\end{equation}
for $\al+i\zeta\in B_{r_0'}$. We take $0<r<r_0'$ and $R_0<R$ to define the open set $O$ as in \eqref{openO}. Therefore we can use the classical method of successive approximations:
$$
(z^{n+1}(t),\varpi^{n+1}(t))=(z^0,\varpi^0)+\int_0^t F(z^n(s),\varpi^n(s))ds,
$$
for $F:O\rightarrow X_{r'}$ and $0\leq r'<r$. We assume by induction that $$\|z^k\|_r(t)< R, \qquad \|\va^k\|_r(t)< R\qquad\mbox{ and }\qquad G(z^k)(\al+i\zeta,\beta,t)> R^{-2}$$ with $\al+i\zeta\in B_r$, $\beta\in\T$ for $k\leq n$  and $0<t<T$ with $T=\min(T_A,T_{CK})$. Here $T_{CK}$ is the time obtaining in the proofs  in \cite{Nirenberg} and \cite{Nishida} (see also \cite{MP}). Now, we will check that $G(z^{n+1})(\al+i\zeta,\beta,t)> R^{-2}$ for $\al+i\zeta\in B_{r}$ and $\beta\in\T$ giving $T_A$. The rest of the proof follows in the same way as in \cite{Nirenberg}, \cite{Nishida}. The following formula:
$$
z^{n+1}(t)=z^0+\int_0^t F_1(z^n(s),\varpi^n(s))ds
$$
yields
\begin{align*}
\begin{split}
G(z^{n+1})(\al+i\zeta,\beta,t)&\geq G(z^0)(\al+i\zeta,\beta)-I_1-2I_2,
\end{split}
\end{align*}
for
$$
I_1=\int_0^t\Big|\frac{F_1(z^n,\varpi^n)(\al+i\zeta,s)-F_1(z^n,\varpi^n)(\al+i\zeta-\beta,s)}{\beta}\Big|^2_* ds
$$
and
$$
I_2=\Big|\frac{z^0(\al\!+\!i\zeta)\!-\!z^0(\al\!+\!i\zeta\!-\!\beta)}{\beta}\Big|_*
\int_0^t\Big|\frac{F_1(z^n,\varpi^n)(\al\!+\!i\zeta,s)\!-\!F_1(z^n,\varpi^n)(\al\!+\!i\zeta\!-\!\beta,s)}{\beta}\Big|_* ds.
$$
Using the induction hypothesis and \eqref{paraarc-chord} it is straightforward to get $I_1\leq C_R^2t$. The inequality
$$
\Big|\frac{z^0(\al\!+\!i\zeta)\!-\!z^0(\al\!+\!i\zeta\!-\!\beta)}{\beta}\Big|_*\leq \sup_{B_{r}}|\da z^0(\al\!+\!i\zeta)|_*< R_0
$$
yields $I_2\leq R_0C_Rt$. Therefore, taking $0<T_A< (R_0^{-2}-R^{-2})(C_R^2+2R_0C_R)^{-1}$, we obtain
$G(z^{n+1})(\al+i\zeta,\beta,t)> R^{-2}$.\\

Proof of Proposition \ref{Lip}: We will show first \eqref{paraarc-chord}. We split as follows
$$
F_1(z,\varpi)(\al+i\zeta)-F_1(z,\varpi)(\al+i\zeta-\beta)=I_1+I_2+I_3
$$
for $$I_1=BR(z,\varpi)(\al+i\zeta)-BR(z,\varpi)(\al+i\zeta-\beta),$$
$$I_2=(H(\varpi)(\al+i\zeta)-H(\varpi)(\al+i\zeta-\beta))\da z(\al+i\zeta),$$
and
$$I_3=H(\varpi)(\al+i\zeta)(\da z(\al+i\zeta)-\da z(\al+i\zeta-\beta)).$$
It is easy to get
$$\sup_{\al+i\zeta\in B_r,\beta\in\T}|I_2|_*\leq \sup_{B_r}|\da z(\al+i\zeta)|_*\sup_{B_r}|H(\da \varpi)(\al+i\zeta)|_*|\beta|,$$
and due to \begin{equation}\label{hcdencd}H:C^{\delta}\rightarrow C^{\delta},\end{equation} (see \cite{St1}), yields
$$\sup_{\al+i\zeta\in B_r,\beta\in\T}|I_2|_* \leq R^2 |\beta|.$$
In a similar fashion it follows:
$$\sup_{\al+i\zeta\in B_r,\beta\in\T}|I_3|_* \leq R^2 |\beta|.$$
For $I_1$, a straightforward calculation gives
$$\sup_{\al+i\zeta\in B_r,\beta\in\T}|I_1|_*\leq \sup_{B_r}|\da BR(z,\varpi)(\al+i\zeta)|_*|\beta|,$$
and it remains to bound $\da BR(z,\varpi)(\al+i\zeta)$. For this term we use the following decomposition:
$$
\da BR(z,\varpi)(\al+i\zeta)=J_1+J_2+J_3,
$$
with
$$
J_1=\frac{1}{2\pi}PV\int_{-\pi}^\pi\da\varpi(\g-\beta)\frac{(z(\g)-z(\g-\beta))^{\bot}}{|z(\g)-z(\g-\beta)|^2}d\beta,
$$
$$
J_2=\frac{1}{2\pi}PV\int_{-\pi}^\pi\varpi(\g-\beta)\frac{\da z(\g)-\da z(\g-\beta)}{|z(\g)-z(\g-\beta)|^2}d\beta,
$$
$$
J_3=-\frac{1}{\pi}PV\int_{-\pi}^\pi\varpi(\g\!-\!\beta)\frac{(z(\g)\!-\!z(\g\!-\!\beta))^\bot}{|z(\g)\!-\!z(\g\!-\!\beta)|^4}
\big((z(\g)\!-\!z(\g\!-\!\beta))\cdot(\da z(\g)\!-\!\da z(\g\!-\!\beta))\big)d\beta,
$$ where $\g=\al+i\zeta$. Here we have to deal with nonlinear singular integral operators given by one-to-one curves. We proceed as in \cite{DY} considering the arc-chord condition (see also \cite{ADY}). we take $J_1=K_1+K_2+K_3$
for
$$
K_1=\frac{1}{2\pi}\int_{-\pi}^\pi\da\varpi(\g-\beta)\frac{(z(\g)-z(\g-\beta)-\da z(\g)\beta)^{\bot}}{|z(\g)-z(\g-\beta)|^2}d\beta,
$$
$$
K_2=\frac{(\da z(\g))^{\bot}}{2\pi}\int_{-\pi}^\pi\da\varpi(\g-\beta)\Big(\frac{\beta}{|z(\g)-z(\g-\beta)|^2}-\frac{1}{|\da z(\g)|^2\beta}\Big)d\beta,
$$
$$
K_3=\frac{(\da z(\g))^{\bot}}{|\da z(\g)|^2}\Big(\frac1{2\pi}\int_{-\pi}^\pi\da\varpi(\g-\beta)[\frac{1}{\beta}-\frac{1}{2\tan(\beta/2)}]d\beta+ H(\da\varpi)(\g) \Big).
$$
We rewrite $K_1$ as follows:
$$
K_1=\frac{1}{2\pi}\int_{-\pi}^\pi\da\varpi(\g-\beta)\frac{(z(\g)-z(\g-\beta)-\da z(\g)\beta)^{\bot}}{
\beta^2}\frac{\beta^2}{|z(\g)-z(\g-\beta)|^2}d\beta,
$$
and therefore, using that $z,\varpi\in O$ and the following estimate:
$$
\sup_{\g\in B_r,\beta\in \T}|z(\g)-z(\g-\beta)-\da z(\g)\beta|_* \leq \sup_{\g\in B_r}|\da^2 z(\g)|_*|\beta|^2,
$$
we obtain $\sup_{B_r}|K_1|_*\leq R^4$. In the integral in $K_2$ we find
$$
\da\varpi(\g-\beta)\Big(\frac{(\da z(\g)\beta+z(\g)-z(\g-\beta))\cdot(\da z(\g)\beta-(z(\g)-z(\g-\beta)))}{|z(\g)-z(\g-\beta)|^2|\da z(\g)|^2\beta}\Big).
$$
The bound for the infimum of $G$ also gives $\sup_{B_r}||\da z(\g)|^{-2}|_*\leq R^2$ in $O$, so we have
$\sup_{B_r}|K_2|_*\leq 2R^8$. The integral in $K_3$ has a bounded kernel in $\beta$ and therefore $$\sup_{B_r}|K_3|_*\leq (C+1) R^4$$ for $C=\max_{\beta\in\T} |\beta^{-1}-(2\tan(\beta/2))^{-1}|$. In $J_2$ we write $J_2=K_4+K_5+K_6$
$$
K_4=\frac{1}{2\pi}\int_{-\pi}^\pi(\varpi(\g-\beta)-\va(\g))\frac{(\da z(\g)-\da z(\g-\beta))^{\bot}}{|z(\g)-z(\g-\beta)|^2}d\beta,
$$
$$
K_5=\frac{\va (\g)}{2\pi}\int_{-\pi}^\pi(\da z(\g)-\da z(\g-\beta))^{\bot} \Big(\frac{1}{|z(\g)-z(\g-\beta)|^2}-\frac{1}{|\da z(\g)|^2\beta^2}\Big)d\beta,
$$
$$
K_6=\frac{\va (\g)}{2|\da z(\g)|^2}\Big(\frac1{\pi}\int_{-\pi}^\pi(\da z(\g)-\da z(\g-\beta))^{\bot}[\frac{1}{\beta^2}-\frac{1}{(2\sin(\beta/2))^2}]d\beta+ (\la(\da z))^{\bot}(\g) \Big),
$$
where $\la=H(\da)$. In $K_4$ we rewrite
$$
K_4=\frac{1}{2\pi}\int_{-\pi}^\pi\frac{(\varpi(\g-\beta)-\va(\g))}{\beta}\frac{(\da z(\g)-\da z(\g-\beta))^{\bot}}{\beta}\frac{\beta^2}{|z(\g)-z(\g-\beta)|^2}d\beta,
$$
and therefore $\sup_{B_r}|K_4|_*\leq R^4$. We take
$$
K_5=\frac{\va (\g)}{2\pi}\int_{-\pi}^\pi\frac{(\da z(\g)-\da z(\g-\beta))^{\bot}}{\beta} \Big(\frac{\beta}{|z(\g)-z(\g-\beta)|^2}-\frac{1}{|\da z(\g)|^2\beta}\Big)d\beta,
$$ and therefore, as for $K_2$, we find $\sup_{B_r}|K_5|_*\leq 2R^8$. The first term in $K_6$ is easy to deal with because the function $\beta^{-2}-(2\sin(\beta/2))^{-2}$ is bounded. For the second term we find
$$
\sup_{B_r}|\la(\da z))^{\bot}|_*=\sup_{B_r}|H(\da^2 z))^{\bot}|_*\leq CR
$$
using \eqref{hcdencd}. For the term $J_3$ we proceed as before to get finally \eqref{paraarc-chord}.

Now we will show how to obtain \eqref{casiL}. The estimate \eqref{cota} follows in a easier fashion (see also \cite{ADY}). Here we will use the following estimate:
\begin{equation}\label{comobm}
\|fg\|_{C^\delta}\leq \|f\|_{C^\delta}\|g\|_{L^\infty}+\|f\|_{L^\infty}\|g\|_{C^\delta}.
\end{equation}

For $F_2$ we get $\da^2F_2(\va)=\da^3(\va H\va)$. In the subtraction $\da^2F_2(\va^2)-\da^2F_2(\va^1)$ we find terms of different order. Here we deal with the singular ones. The rest of the terms can be estimate in a simpler way. We find in the subtraction the term
$I_1=\da^3\va^2 H(\va^2)-\da^3\va^1 H(\va^1)$ and we split
$$
I_1=(\da^3\va^2-\da^3\va^1)H(\va^2)+\da^3\va^1(H(\va^2)-H(\va^1))=J_1(\g)+J_2(\g).
$$
For $\beta\in \T$ and $\g\in B_{r'}$, the inequality \eqref{comobm} yields
\begin{align*}
\frac{|J_1(\g)-J_1(\g-\beta)|_*}{|\beta|^{\delta}}&\leq 2\|\da\va^2-\da\va^1\|_{r'}\|\va^2\|_{r'},
\end{align*}
and using \eqref{Cauchy} it follows:
$$
\sup_{\g\in B_{r'},\beta\in\T}\frac{|J_1(\g)-J_1(\g-\beta)|_*}{|\beta|^{\delta}}\leq \frac{2R}{r-r'}\|\va^2-\va^1\|_{r}.
$$
Analogously
$$
\frac{|J_2(\g)\!-\!J_2(\g\!-\!\beta)|_*}{|\beta|^{\delta}}\leq 2 \|\da\va^1\|_{r'}\|\va^2\!-\!\va^1\|_{r'} \leq \frac{2\|\va^1\|_{r}}{r-r'}\|\va^2\!-\!\va^1\|_{r'}\leq \frac{2R}{r-r'}\|\va^2\!-\!\va^1\|_{r}.
$$
Also the term $I_2$, given by $I_2=\va^2 H(\da^3\va^2)-\va^1 H(\da^3\va^1)$, can be decomposed as
$$
J_3(\g)=(\va^2-\va^1)H(\da^3\va^2),\quad J_4(\g)=\va^1H(\da^3(\va^2-\va^1),
$$
and as before
\begin{align*}
\frac{|J_3(\g)-J_3(\g-\beta)|_*}{|\beta|^{\delta}}&\leq 2\|\va^2-\va^1\|_{r'}\|H(\da\va^2)\|_{r'}\leq 2C\|\va^2-\va^1\|_{r'}\|\da\va^2\|_{r'}\\
&\leq \frac{2RC}{r-r'}\|\va^2-\va^1\|_{r}.
\end{align*}
For $J_4$ it follows:
\begin{align*}
\frac{|J_4(\g)-J_4(\g-\beta)|_*}{|\beta|^{\delta}}&\leq 2RC\|\da(\va^2-\va^1)\|_{r'}\leq \frac{2RC}{r-r'}\|\va^2-\va^1\|_{r}.
\end{align*}
Now we consider the operator $F_1(z,\va)=BR(z,\va)+H\va \da z$. The estimates for the second term are as before, we then show the control for the Birkhoff-Rott integral. While the terms of lower order are easier, we consider in $\da^2BR(z,\va)$ the most  singular:
$$
I_3=\frac{1}{2\pi}PV\int_{-\pi}^\pi\varpi(\g-\beta)\frac{(\da^2 z(\g)-\da^2 z(\g-\beta))^{\bot}}{|z(\g)-z(\g-\beta)|^2}d\beta,
$$
$$
I_4=\frac{-1}{\pi}PV\int_{-\pi}^\pi\varpi(\g\!-\!\beta)\frac{(z(\g)\!-\!z(\g\!-\!\beta))^{\bot}}{|z(\g)\!-\!z(\g\!-\!\beta)|^4}
(z(\g)\!-\!z(\g\!-\!\beta))\cdot(\da^2 z(\g)\!-\!\da^2 z(\g\!-\!\beta))d\beta,
$$
and
$$
I_5=\frac{1}{2\pi}PV\int_{-\pi}^\pi\da^2\varpi(\g-\beta)\frac{(z(\g)-z(\g-\beta))^{\bot}}{|z(\g)-z(\g-\beta)|^2}d\beta.
$$
We take $I_3=J_5+J_6+J_7+J_8$ with
$$
J_5=\frac{1}{2\pi}\int_{-\pi}^\pi(\varpi(\g-\beta)-\va(\g))\frac{(\da^2 z(\g)-\da^2 z(\g-\beta))^{\bot}}{|z(\g)-z(\g-\beta)|^2}d\beta,
$$
$$
J_6=\frac{\va(\g)}{2\pi}\int_{-\pi}^\pi(\da^2 z(\g)-\da^2 z(\g-\beta))^{\bot}[\frac{1}{|z(\g)-z(\g-\beta)|^2}-\frac{1}{|\da z(\g)|^2\beta^2}]d\beta,
$$
$$
J_7=\frac{\va(\g)}{2\pi|\da z(\g)|^2}\int_{-\pi}^\pi(\da^2 z(\g)-\da^2 z(\g-\beta))^{\bot}[\frac{1}{\beta^2}-\frac{1}{4\sin^2(\beta/4)}]d\beta,
$$
and
$$
J_8=\frac{\va(\g)}{2|\da z(\g)|^2}(\la(\da^2 z)(\g))^{\bot}.
$$
Then, with this splitting, in $\da^2BR(z^2,\va^2)-\da^2BR(z^1,\va^1)$, one can find the term
$$
DJ_8=\frac{\va^2(\g)}{2|\da z^2(\g)|^2}(\la(\da^2 z^2)(\g))^{\bot}-\frac{\va^1(\g)}{2|\da z^1(\g)|^2}(\la(\da^2 z^1)(\g))^{\bot}
$$
Now, for $h\in \T$ and $\g\in B_{r'}$, it follows:
\begin{align*}
|DJ_8(\g)-DJ_8(\g-h)|_*&\leq C_R(\|(z^2-z^1,\va^2-\va^1)\|_{r'}|h|^\delta(r-r')^{-1}\\
&\quad\quad\quad +|\la(\da^2(z^2-z^1))(\g)-\la(\da^2(z^2-z^1))(\g-h)|_*),
\end{align*}
and using \eqref{hcdencd} one finds
\begin{align*}
|\la(\da^2(z^2\!-\!z^1))(\g)\!-\!\la(\da^2(z^2\!-\!z^1))(\g\!-\!h)|_*&=|H(\da^3(z^2\!-\!z^1))(\g)\!-\!H(\da^3(z^2\!-\!z^1))(\g\!-\!h)|_*\\
&\leq C\|\da(z^2-z^1)\|_{r'}|h|^\delta,
\end{align*}
and finally
$$
\frac{|DJ_8(\g)-DJ_8(\g-h)|_*}{|h|^\delta}\leq \frac{C_R}{r-r'}\|(z^2-z^1,\va^2-\va^1)\|_{r}.
$$
In an analogous way we may define $DJ_5$ and split it as follows:
$$
K_7=\frac{1}{2\pi}\int_{-\pi}^\pi((\va^2-\va^1)(\g-\beta)-(\va^2-\va^1)(\g))\frac{(\da^2 z^2(\g)-\da^2 z^2(\g-\beta))^{\bot}}{|z^2(\g)-z^2(\g-\beta)|^2}d\beta,
$$
$$
K_8=\frac{1}{2\pi}\int_{-\pi}^\pi(\va^1(\g-\beta)-\va^1(\g))\frac{(\da^2 (z^2-z^1)(\g)-\da^2 (z^2-z^1)(\g-\beta))^{\bot}}{|z^2(\g)-z^2(\g-\beta)|^2}d\beta,
$$
and
$$
K_9=\frac{1}{2\pi}\int_{-\pi}^\pi(\va^1(\g-\beta)-\va^1(\g))(\da^2 z^1(\g)-\da^2 z^1(\g-\beta))^{\bot}A(\g,\beta)d\beta,
$$
with $A(\g,\beta)=|z^2(\g)-z^2(\g-\beta)|^{-2}-|z^1(\g)-z^1(\g-\beta)|^{-2}$. All kernels in the integrals in $K_7$, $K_8$ and $K_9$ have grade $0$ so the control of all these terms are analogous. Now we will show in detail the term $K_7$. We rewrite it as
$$
K_7=\frac{1}{2\pi}\int_{-\pi}^\pi B(\g,\beta)C(\g,\beta) D(\g,\beta)d\beta,
$$
with
$$
B(\g,\beta)=\frac{(\va^2\!-\!\va^1)(\g\!-\!\beta)\!-\!(\va^2\!-\!\va^1)(\g)}{\beta}$$
$$
C(\g,\beta)=\frac{(\da^2 z^2(\g)\!-\!\da^2 z^2(\g\!-\!\beta))^{\bot}}{\beta},\qquad D(\g,\beta)=\frac{\beta^2}{|z^2(\g)\!-\!z^2(\g\!-\!\beta)|^2},
$$
to get the following splitting
$$
K_7(\g)-K_7(\g-h)=L_1+L_2+L_3,
$$
where
$$
L_1=\frac{1}{2\pi}\int_{-\pi}^\pi (B(\g,\beta)-B(\g-h,\beta))C(\g,\beta) D(\g,\beta)d\beta,
$$
$$
L_2=\frac{1}{2\pi}\int_{-\pi}^\pi B(\g,\beta)(C(\g,\beta)-C(\g-h,\beta))D(\g,\beta)d\beta,
$$
and
$$
L_3=\frac{1}{2\pi}\int_{-\pi}^\pi B(\g,\beta)C(\g,\beta)(D(\g,\beta)-D(\g-h,\beta))d\beta.
$$
We take the term $B$ as
$$
B(\g,\beta)=\int_0^1\da(\va^2-\va^1)(\g-s\beta)\,ds,
$$
and therefore
$$
|B(\g,\beta)-B(\g-h,\beta)|_*\leq \|\va^2-\va^1\|_{r'}|h|^\delta.
$$
for $\g\in B_{r'}$ and $h\in\T$. It yields
$$
|L_1|_*\leq \|\va^2-\va^1\|_{r'}|h|^\delta\|\da z^2\|_{r'}R^2\leq
\frac{R^3}{r-r'} \|\va^2-\va^1\|_{r}|h|^\delta.
$$
For $C(\g,\beta)$ it follows:
$$
C(\g,\beta)=\int_0^1(\da^3 z^2(\g+(s-1)\beta))^{\bot}\,ds,
$$
and analogously one gets
$$
|L_2|_*\leq \|\va^2-\va^1\|_{r'}\|\da z^2\|_{r'}|h|^\delta R^2\leq \frac{R^3}{r-r'} \|\va^2-\va^1\|_{r}|h|^\delta.
$$
In $L_3$ we rewrite the difference $D(\g,\beta)-D(\g-h,\beta)$ as
\begin{align*}
\frac{\beta^2}{|z^2(\g)\!-\!z^2(\g\!-\!\beta)|^2}
\frac{\beta^2}{|z^2(\g-h)\!-\!z^2(\g-h\!-\!\beta)|^2}
E_1(\g,h,\beta)\cdot E_2(\g,h,\beta),
\end{align*}
where
$$
E_1(\g,h,\beta)=\frac{(z^2(\g\!-\!h)\!-\!z^2(\g\!-\!h\!-\!\beta))\!+\!(z^2(\g)\!-\!z^2(\g\!-\!\beta))}{\beta}
$$
$$
E_2(\g,h,\beta)=\frac{(z^2(\g\!-\!h)\!-\!z^2(\g\!-\!h\!-\!\beta))\!-\!(z^2(\g)\!-\!z^2(\g\!-\!\beta))}{\beta}
$$
As before one can take
$$
E_2(\g,h,\beta)=\int_0^1\big(\da z^2(\g\!-\!h+(s-1)\beta)-\da z^2(\g\!+(s-1)\beta)\big)
$$
and therefore $|E_2|_*\leq \|z^2\|_{r'}|h|^\delta$. It provides as before
$$
|D(\g,\beta)-D(\g-h,\beta)|_*\leq 2R^6|h|^\delta,
$$
and
$$
|L_3|_*\leq \|\va^2-\va^1\|_{r'}\|\da z^2\|_{r'}2R^6|h|^\delta \leq \frac{ 2R^7}{r-r'}\|\va^2-\va^1\|_{r}|h|^\delta.
$$
All these estimates for the terms $L_j$ yield
$$
\frac{|K_7(\g)-K_7(\g-h)|_*}{|h|^\delta}\leq \frac{C_R}{r-r'}\|\va^2-\va^1\|_{r},
$$
for $\g\in B_r'$ and $h\in\T$.

In a similar way it is possible to get the appropriate control for $J_6$ and $J_7$. The terms $I_4$ and $I_5$ can be estimated as $I_3$, so that with this argument we finish the proof.

\section{Ill-posedness for the amplitude equation.}

In this section we choose the tangential term  $c(\alpha,t)=\frac{1}{2}H\va(\alpha,t)$ which gives the following closed equation for the amplitude of the vorticity
\begin{eqnarray}
\va_t-\frac{1}{2}(\va H\va)_\si=0,\label{oneD}\\
\va(\si,0)=\va_0(\si).
\end{eqnarray}
  We shall prove the following theorem:

\begin{thm} Let  $\va_0\in H^s(\ff{T})$  with $s>\frac{3}{2}$ and $$\int_\ff{T}\va_0=0.$$

Then if there exists a point $\si_0$ where $\va_0(\si_0)>0$ and $\va_0$ is not $C^\infty$ in $\si_0$, there is no solution of equation (\ref{oneD}) in the class $C([0,T);H^s(\ff{T}))$ with $s>\frac{3}{2}$ and $T>0$. In addition, $\va_0\in C^\infty$ is not sufficient to obtain existence.
\end{thm}
\begin{rem}
  In the case of the real line $\ff{R}$ equation (\ref{oneD}) is also ill-posed, in $H^s$ with $s>3/2$, for a non-analytic initial data. For more details see \cite{AD}.
\end{rem}
Proof: We will proceed by a contradiction argument.

Let us assume that there exist a solution of equation
(\ref{oneD}) in the class $C([0,T),H^{s}(\mathbb{T}))$ with $\va(\si,0)=\va_0(\si)$.

First we have to note that if the initial data $\va_0$ are of mean zero, then the solution $\va$ will remain of mean zero.

Now, taking the Hilbert transform on equation (\ref{oneD})  yields
\begin{equation*}
\pa_tH\va-\frac{1}{2}(H\va H\va_\si-\va \va_\si)=0,
\end{equation*}
where we have used the following properties of the Hilbert transform for a periodic function with mean zero:
\begin{itemize}
\item $H(H\va)=-\va.$
\item $H(\va H\va)=\frac{1}{2}((H\va)^2-\va^2).$
\end{itemize}
 We denote the complex valued function
$z(\si,t)=H\va (\si,t)-i\va (\si,t)$  which satisfies
\begin{equation}
\pa_tz-\frac{1}{2}zz_\si=0\label{burger}.
\end{equation}
  Take $P_{\si}(u)$ to be the Green's function of the Laplacian  for the Dirichlet problem in the unit ball
$$P_{\si}(u)\equiv\frac{1}{2\pi}\frac{1-|u|^2}{|u-\si|^2},$$
and $P\va (u)$ will be
$$P\va (u)\equiv \int_{\pa B(0,1)}P_{\tau}(u)\va(\tau) d\tau.$$
Therefore
$$Z(u)=P(H\va-i\va)(u) \quad\text{with} \quad u=r e^{i\si},$$
is an analytic function on the unit ball.
Applying $P$ to the equation (\ref{burger}) yields
$$\pa_tPz=\frac{1}{2}P(zz_\si),$$
where we can write the second term in the following way
$$P(zz_\si)=Pz(Pz)_\si,$$
since both terms have the same restriction to the boundary of the unit ball and both are harmonic.

Thus, we have for $Z(u,t)$  the equation
$$Z_t-\frac{1}{2}ZZ_\si=0\quad \text{on}\quad u\in \overline{B(0,1)},$$
hence
\begin{eqnarray}
Z_t-\frac{1}{2}iuZZ_u=0\quad \text{on} \quad u\in B(0,1),\label{complexburger}\\
Z(u,0)=Z_0(u)=P(H\va_0-i\va_0)(u).
\end{eqnarray}
We will define the complex trajectories $X(u,t)$ by
$$\frac{dX(u,t)}{dt}=-\frac{1}{2}iX(u,t)Z(X(u,t),t),$$
$$X(u,0)=u,\quad u\in B(0,1).$$
For
sufficiently small $t$, by Picard's Theorem, these trajectories
exist and $X(u,t)\in B(0,1).$ Therefore
$$\frac{d Z(X(u,t),t)}{dt}=\pa_tZ(X(u,t),t)-\frac{1}{2}iX(u,t)Z(X(u,t),t)Z_u(X(u,t),t)=0.$$
Thus, we have
$$Z(X(u,t),t)=Z_0(u),$$
and
$$\frac{dX(u,t)}{dt}=-\frac{1}{2}iX(u,t)Z_0(u).$$
Moreover
$$X(u,t)=ue^{-\frac{1}{2}iZ_0(u)t}.$$
Taking modules in the last expression we obtain
$$R(u,t)=|X(u,t)|=re^{-\frac{1}{2}P\va_0(re^{i\si})t}.$$
If we consider a point $e^{i\sigma_0}=u_0\in \pa B(0,1)$ with $\va_0(\si_0)>0$, then
$$R(u_0,t)=e^{-\frac{1}{2}w_0(\si_0)t}<1.$$
Hence $X(u_0,t)\in B(0,1)$ for all $t>0$, and a continuity argument yields
$$Z(X(\si_0,t),t)=z_0(\si_0)=H\va_{0}(\si_0)-i\va_0(\si_0),$$
where to simplify we denote $X(u_0,t)=X(\si_0,t)$. Then we have
$$X(\si_0,t)=e^{i(\si_0-\frac{1}{2}z_0(\si_0)t)}.$$
Taking a derivative with respect to $\si_0$ on this equation we find that
$$\frac{dX(\si_0,t)}{d\si_0}=i(1-\frac{1}{2}z_{0\,\si}(\si_0)t)X(\si_0,t).$$
With the chain's rule we obtain
$$\frac{dZ}{dX}(X(\si_0,t),t)iX(\si_0,t)=\frac{dZ}{d\Theta}(X(\si_0,t),t)=\frac{z_{0\,\si}(\si_0)}{(1-\frac{1}{2}z_{0\,\si}(\si_0)t)},$$
where
$$X(\si_0,t)=R(\si_0,t)e^{i\Theta(\si_0,t)}.$$
Taking two derivatives
$$\frac{d^2Z}{d\Theta^2}(X(\si_0,t),t)=\frac{z_{0\,\si\si}(\si_0)}{(1-\frac{1}{2}z_{0\,\si}(\si_0)t)^3}.$$
For the n-th derivative we have
$$\frac{d^nZ}{d\Theta^n}(X(\si_0,t),t)=\frac{\frac{d^nz_0}{d\si^n}(\si_0)}{(1-\frac{1}{2}z_{0\,\si}(\si_0)t)^{n+1}}+\text{``lower terms''}.$$
We observe that $(1-\frac{1}{2}z_{0\,\si}(\si_0)t)\neq 0$ for $t$ small enough.

Then if $w_0$ is not $C^\infty$ in $\si_0$ this is a contradiction since $Z(u,t)$ is analytic on $X(\si_0,t)$ for all $t>0$.

In addition, if
$\va_0(\si_0)>0$ and $\frac{d^n\va_0}{d\si^n}(\si_0)=0\,\,\forall n$ but $\va_0$ is
not constant on any neighborhood of $\si_0$, we can conclude
$$\frac{d \Im Z}{d\Theta}(X^1(\si_0,t),X^2(\si_0,t))=0.$$
Continuing this process we obtain that all derivatives satisfy
$$\frac{d^n\Im Z}{d\Theta^n}(X^1(\si_0,t),X^2(\si_0,t))=0.$$
 The imaginary part $\Im Z(x_1,x_2,t)$ is analytic on $(x_1,x_2)=(X^1(\si_0,t),X^2(\si_0,t))$ for all $t>0$, thus $\Im Z(x_1,x_2)$ is constant over the circumference, $R=R(\si_0,t)$, and this is a contradiction if $\va_0$ is not constant .


\section{Appendix}


Here we extend the property of the continuity of the pressure
known for Darcy's flow (see \cite{DY} and \cite{ADY}). Writing
this work we learned of the paper by Shvydkoy \cite{Roman} who
also obtain this fact for more general cases in a different way.

\begin{prop}\label{presiones}
Let us consider a weak solution $(v,p)$ satisfying
(\ref{weakeuler}--\ref{weakdiv0}) where $\curl v=\omega$ is given
by \eqref{omegacomoes}. Then we have the following identity
$$p^1(z(\al,t),t)=p^2(z(\al,t),t),$$ where $p^j(z(\alpha,t),t)$
denotes the limit pressure obtained approaching the boundary in
the normal direction inside $\Omega^j$.
\end{prop}

Proof: We shall show that the Laplacian of the pressure is as
follows
$$
\Delta p(x,t)= F(x,t)+f(\alpha,t)\delta(x-z(\al,t)),
$$
where $F$ is regular in $\Omega^j(t)$ and discontinuous on
$z(\al,t)$. The amplitude of the delta function $f$ is regular.
The inverse of the Laplacian by means of the Newtonian potential
gives the continuity of the pressure on the free boundary (see
\cite{ADY}). Here we shall give the argument for a close curve;
the proof for the other cases being analogous.

The expression for the conjugate of the velocity in complex
variables
$$
\overline{v}(z,t)=\frac{1}{2\pi i}PV\int
\frac{1}{z-z(\al,t)}\varpi(\al,t)d\al,
$$
for $z\neq z(\al,t)$ allows us to accomplish the fact that
\begin{align*}
\partial_z\overline{v}(z,t)&=\frac{1}{2\pi i}PV\int \frac{-\varpi(\al,t)}{(z-z(\al,t))^2}d\al=\frac{1}{2\pi i}PV\int \frac{-\da z(\al,t)}{(z-z(\al,t))^2}\frac{\varpi(\al,t)}{\da z(\al,t)}d\al.
\end{align*}
Therefore
\begin{equation}\label{dzv}
\partial_z\overline{v}(z,t)=\frac{1}{2\pi i}PV\int \frac{1}{z-z(\al,t)}\da(\frac{\varpi}{\da z})(\al,t) d\al
\end{equation}
for a regular parametrization with $\da z(\al,t)\neq 0$. In a
similar way
\begin{align}\label{d2zv}
\partial^2_z\overline{v}(z,t)&=\frac{1}{2\pi i}PV\int \frac{1}{z-z(\al,t)}\da(\frac{1}{\da z}\da(\frac{\varpi}{\da z}))(\al,t) d\al.
\end{align}
These identities allow us to get the values of $\grad
v^j(z(\al,t),t)$ and $\grad^2 v^j(z(\al,t),t)$. As for the
velocity, the limits are different, but we can compute the values.

To get the above formula for the pressure we take the weak type
identity \eqref{weakeuler} with $\eta(x,t)=\grad\lambda(x,t)$. We
can compute then the Laplacian of the pressure in a weak sense
due to
\begin{align*}
\int_0^T\!\!\int_{\R^2}p\Delta\lambda dxdt&=-
\int_0^T\!\!\int_{\R^2} v\cdot \grad\lambda_t dxdt-\int_0^T\!\!\int_{\R^2}  v\cdot(v\cdot\grad^2\lambda)dxdt-\int_{\R^2}v_0(x)\cdot\grad\lambda(x,0)dx\\
&=I_1+I_2+I_3.
\end{align*}
Then
$$I_1=I_3=0$$
by the incompressible condition. We define
$$\Omega_{\ep}^1(t)=\{x\in\Omega^1(t):\dist(x,\partial\Omega^1(t))\geq\ep\}$$
$$\Omega_{\ep}^2(t)=\{x\in\Omega^2(t):\dist(x,\partial\Omega^2(t))\geq\ep\}.$$
We decompose as follows $I_2=J_3+J_4+J_5+J_6$ as previously where
$$
J_3=-\int_0^T\!\!\int_{\R^2}(v_1)^2\dxu^2\lambda dxdt,\qquad
J_4=-\int_0^T\!\!\int_{\R^2} v_1v_2\dxd\dxu\lambda dxdt,
$$
$$
J_5=-\int_0^T\!\!\int_{\R^2} v_1v_2\dxu\dxd\lambda dxdt,\qquad
J_6=-\int_0^T\!\!\int_{\R^2}(v_2)^2\dxd^2\lambda dxdt.
$$
Using the sets $\Omega^j_\ep(t)$ and the identity \eqref{dzv} we
get
\begin{align*}
J_3&=-\lim_{\ep\rightarrow 0}(\int_0^T\!\!\int_{\Omega^1_{\ep}(t)}(v_1)^2\dxu^2\lambda dxdt+\int_0^T\!\!\int_{\Omega^2_{\ep}(t)}(v_1)^2\dxu^2\lambda dxdt)\\
&=\int_0^T\!\!\int_{\R^2}2 v_1\dxu v_1\dxu\lambda dxdt\\
&\quad+\int_0^T\!\!\int_{-\pi}^\pi\big((v^2_1(z(\al,t),t))^2-(v_1^1(z(\al,t),t))^2\big)\dxu\lambda(z(\al,t),t)\da z_2(\al,t) d\al dt\\
&=K_1+K_2.
\end{align*}
The term $K_1$ trivializes because the subtle integration by parts
and the identity \eqref{d2zv} give
\begin{align*}
K_1&=-\int_0^T\!\!\int_{\R^2}2 (v_1\dxu^2 v_1+(\dxu v_1)^2)\lambda
dxdt-\int_0^T\!\!\int_{-\pi}^\pi
\widetilde{f}(\al,t)\lambda(z(\al,t),t) d\al dt
\end{align*}
for $\widetilde{f}(\al,t)=2(v^2_1(z(\al,t),t)\dxu
v^2_1(z(\al,t),t)- v_1^1(z(\al,t),t) \dxu v^1_1(z(\al,t),t))\da
z_2(\al,t)$. The first term in $K_1$ is part of $F(x,t)$ and the
second of $f(\al,t)$.

We can rewrite $K_2$ as follows
\begin{align}
\begin{split}
K_2&=-2\int_0^T\!\!\int_{-\pi}^\pi \varpi BR_1\frac{\da z_1}{|\da
z|^2}\dxu\lambda(z)\da z_2 d\al dt.
\end{split}
\end{align}
We continue with $J_4$
\begin{align*}
J_4&=\int_0^T\!\!\int_{\R^2} (v_2\dxd v_1 +v_1\dxd v_2)\dxu\lambda dxdt\\
&\quad-\int_0^T\!\!\int_{-\pi}^\pi\big((v^2_1v^2_2)(z(\al,t),t)-(v_1^1v_2^1)(z(\al,t),t)\big)\dxu\lambda(z(\al,t),t)\da z_1(\al,t) d\al dt\\
&=K_3+K_4.
\end{align*}
We deal with the term $K_3$ in a similar way as with $K_1$.

We can rewrite $K_4$ as follows
\begin{align}
\begin{split}
K_4=-\int_0^T\!\!\int_{-\pi}^\pi [\varpi BR_1\frac{\da z_2}{|\da
z|^2}+\varpi BR_2\frac{\da z_1}{|\da z|^2}]\dxu\lambda(z)\da z_1
d\al dt.
\end{split}
\end{align}
For $J_5$ we split
\begin{align*}
J_5&=\int_0^T\!\!\int_{\R^2} (v_2\dxu v_1 +v_1\dxu v_2)\dxd\lambda dxdt\\
&\quad+\int_0^T\!\!\int_{-\pi}^\pi\big((v^2_1v^2_2)(z(\al,t),t)-(v_1^1v_2^1)(z(\al,t),t)\big)\dxd\lambda(z(\al,t),t)\da z_2(\al,t) d\al dt\\
&=K_5+K_6.
\end{align*}
We proceed for $K_5$ in a similar manner  as with $K_1$.

We obtain for $K_6$ the following expression
\begin{align}
\begin{split}
K_6&=\int_0^T\!\!\int_{-\pi}^\pi [\varpi BR_1\frac{\da z_2}{|\da
z|^2}+\varpi BR_2\frac{\da z_1}{|\da z|^2}]\dxd\lambda(z)\da z_2
d\al dt.
\end{split}
\end{align}
With $J_6$ one finds
\begin{align*}
J_6&=\int_0^T\!\!\int_{\R^2}2 v_2\dxd v_2\dxd\lambda dxdt\\
&\quad-\int_0^T\!\!\int_{-\pi}^\pi\big((v^2_2(z(\al,t),t))^2
-(v_2^1(z(\al,t),t))^2\big)
\dxd\lambda(z(\al,t),t)\da z_2(\al,t) d\al dt\\
&=K_7+K_8.
\end{align*}
For $K_7$ we proceed as before. We obtain for $K_8$ the following
expression
\begin{align}
\begin{split}
K_8&=-2\int_0^T\!\!\int_{-\pi}^\pi \varpi BR_2\frac{\da z_2}{|\da
z|^2}\dxd\lambda(z)\da z_1 d\al dt.
\end{split}
\end{align}
We now sum as follows $K_2+K_4+K_6+K_8=L_2$, then

$$
L_2=-\int_0^T\!\!\int_{-\pi}^\pi \frac{\varpi(\al,t)}{|\da
z(\al,t)|^2} BR(z,\varpi)(\al,t)\cdot\dpa z(\al,t) \da
z(\al,t)\cdot \grad\lambda(z(\al,t),t) d\al dt.
$$
An integration by parts in the variable $\al$ in $L_2$ gives the last term of $f(\al,t)$. The formula for the Laplacian of $p$ is found. \quad\eop\\

\subsection*{{\bf Acknowledgements}}

\smallskip

The authors were partially supported by the grant {\sc MTM2008-03754} of the MCINN (Spain) and
the grant StG-203138CDSIF  of the ERC.

\begin{quote}
\begin{tabular}{l}
\textbf{Angel Castro} \\
{\small Instituto de Ciencias Matem\'aticas}\\
{\small Consejo Superior de Investigaciones Cient\'ificas}\\
{\small Serrano 123, 28006 Madrid, Spain}\\
{\small Email: angel\underline{  }castro@icmat.es}
\end{tabular}
\end{quote}
\begin{quote}
\begin{tabular}{ll}
\textbf{Diego C\'ordoba} &  \textbf{Francisco Gancedo}\\
{\small Instituto de Ciencias Matem\'aticas} & {\small Department of Mathematics}\\
{\small Consejo Superior de Investigaciones Cient\'ificas} & {\small University of Chicago}\\
{\small Serrano 123, 28006 Madrid, Spain} & {\small 5734 University Avenue, Chicago, IL 60637}\\
{\small Email: dcg@icmat.es} & {\small Email: fgancedo@math.uchicago.edu}
\end{tabular}
\end{quote}

\end{document}